\documentclass{amsproc}
\usepackage{amssymb}
\usepackage{amscd}
\usepackage{amsmath}

\theoremstyle{plain}

\newtheorem{remark}{Remark}

\newtheorem{theorem}{Theorem}
\numberwithin{equation}{section}

\def\Nset{\mbox{I\kern-.21em N}}
\def\RE{{\mbox{\rm I\kern-.21em R}}}
\def\ZZ{{\mbox{\sf Z\kern-.45em Z}}}
\def\vv{\kern.344em{\rule[.18ex]{.075em}{1.32ex}}\kern-.344em}

  \def\g{\gamma}  
  
 \def\f{\varphi}

\def\<{\langle} \def\>{\rangle}

\begin{document}

\title[On an inverse dynamic problem for the wave equation]{On an inverse dynamic problem for the wave equation with a potential on a real line.}
\author{ A. S. Mikhaylov}
\address{St. Petersburg   Department   of   V.A. Steklov    Institute   of   Mathematics
of   the   Russian   Academy   of   Sciences, 7, Fontanka, 191023
St. Petersburg, Russia and Saint Petersburg State University,
St.Petersburg State University, 7/9 Universitetskaya nab., St.
Petersburg, 199034 Russia.} \email{mikhaylov@pdmi.ras.ru}
\author{ V. S. Mikhaylov}
\address{St.Petersburg   Department   of   V.A.Steklov    Institute   of   Mathematics
of   the   Russian   Academy   of   Sciences, 7, Fontanka, 191023
St. Petersburg, Russia and Saint Petersburg State University,
St.Petersburg State University, 7/9 Universitetskaya nab., St.
Petersburg, 199034 Russia.} \email{ftvsm78@gmail.com}


\maketitle

\begin{center}
{\bf Abstract.} We consider the inverse dynamic problem for the
wave equation with a potential on a real line. The forward
initial-boundary value problem is set up with a help of boundary
triplets. As an inverse data we use an analog of a response
operator (dynamic Dirichlet-to-Neumann map). We derive equations
of inverse problem and also point out the relationship between
dynamic inverse problem and spectral inverse problem from a
matrix-valued measure.
\end{center}

\section{Introduction}

For a potential $q\in C^2(R)\cap L_1(\mathbb{R})$ we consider an
operator $H$ in $L_2(\mathbb{R})$ given by
\begin{eqnarray*}
(Hf)(x)=-f''(x)+q(x)f(x),\quad x\in \mathbb{R},\\
\operatorname{dom}H=\left\{f\in H^2(\mathbb{R})\,|\,
f(0)=f'(0)=0\right\}.
\end{eqnarray*}
Then
\begin{eqnarray*}
(H^*f)(x)=-f''(x)+q(x)f(x),\quad x\in \mathbb{R},\\
\operatorname{dom}H^*=\left\{f\in L_2(\mathbb{R})\,|\, f\in
H^2(-\infty,0),\,f\in H^2(-\infty,0)\right\}.
\end{eqnarray*}
For a continuous function $g$ we denote
\begin{equation*}
g_\pm:=\lim_{\varepsilon\to 0}g(0\pm\varepsilon).
\end{equation*}
Let $B:=\mathbb{R}^2$. The \emph{boundary operators}
$\Gamma_{0,1}: \operatorname{dom}H^*\mapsto B$ are introduced by
the rules
\begin{equation*}
\Gamma_0w:=\begin{pmatrix}w_+-w_-\\ w'_+-w'_-\end{pmatrix},\quad
\Gamma_1w:=\frac{1}{2}\begin{pmatrix}w'_++w'_-\\
-w_+-w'_-\end{pmatrix}.
\end{equation*}
Integrating by parts for $u,v\in \operatorname{dom}H^*$ shows that
the abstract second Green identity holds:
\begin{equation*}
\left(H^*u,v\right)_{L_2(\mathbb{R})}-\left(u,H^*v\right)_{L_2(\mathbb{R})}=\left(\Gamma_1u,\Gamma_0v\right)_B-\left(\Gamma_0u,\Gamma_1v\right)_B.
\end{equation*}
The mapping
\begin{equation*}
\Gamma:=\begin{pmatrix}\Gamma_0\\\Gamma_1\end{pmatrix}:\operatorname{dom}H^*\mapsto
B\times B
\end{equation*}
evidently is surjective. Then a triplet $\{B,\Gamma_0,\Gamma_1\}$
is a \emph{boundary triplet} for $H^*$ (see \cite{BMN}). With the
help of boundary triplets one can describe self-adjoint extensions
of $H$, see \cite{DM,R,Koch}. In \cite{BD} the authors used the
concept of boundary triplets to set up and study a boundary value
problem for abstract dynamical system with boundary control in
Hilbert space, they also used it for the purpose of describing the
special (wave) model of the one-dimensional Schr\"odinger operator
on an interval \cite{BS}.

Let $T>0$ be fixed. We use the triplet $\{B,\Gamma_0,\Gamma_1\}$
to set up the dynamical system with special boundary control
(acting in the origin) for a wave equation with a potential on a
real line:
\begin{eqnarray}
\label{wave_eqn}u_{tt}+H^*u=0,\quad t> 0,\\
\label{bnd_eqn}(\Gamma_0u)(t)=\begin{pmatrix}f_1(t) \\ f_2'(t)\end{pmatrix}, \quad t>0,\\
\label{init_eqn}u(\cdot,0)=u_t(\cdot,0)=0.
\end{eqnarray}
Here the function $F=\begin{pmatrix}f_1 \\ f_2\end{pmatrix}$,
$f_1\, f_2\in L_2(0,T),$ is interpreted as a \emph{boundary
control}. The solution to (\ref{wave_eqn})--(\ref{init_eqn}) is
denoted by $u^F.$ The \emph{response operator}, the analog of a
Dirichlet-to-Neumann map is introduced by the rule
\begin{equation*}
\left(R^TF\right)(t):=\left(\Gamma_1u^F\right)(t),\quad t>0.
\end{equation*}
The speed of the wave propagation in the system
(\ref{wave_eqn})--(\ref{init_eqn}) equal to one, that is why the
natural set up of the dynamic inverse problem is to find a
potential $q(x),$ $x\in (-T,T)$ from the knowledge of a response
operator $R^{2T}$ (cf. \cite{B07,BM01,AM}).

In the second section we derive the representation formula for the
solution $u^F$ and introduce the operators of the Boundary Control
 method. In the third section we derive Krein and
Gelfand-Levitan equations of the dynamic inverse problem and point
out the the relationship between dynamic and spectral inverse
problems.

\section{Forward problem, operators of the Boundary Control method}

It is straightforward to check that when $q=0$, the solution to
(\ref{wave_eqn})--(\ref{init_eqn}) is given by:
\begin{eqnarray*}
u^F(x,t)=\left\{\begin{array}l
\frac{1}{2}f_1(t-x)-\frac{1}{2}f_2(t-x),\quad x>0,\\
-\frac{1}{2}f_1(t+x)-\frac{1}{2}f_2(t+x),\quad x<0,\\
0, \quad 0<t<|x|.
\end{array}
\right.
\end{eqnarray*}
Everywhere we consider operators acting in $L_2-$spaces, that is
why it is reasonable to introduce the \emph{outer space} of the
system (\ref{wave_eqn})--(\ref{init_eqn}), the space of controls
as $\mathcal{F}^T:= L_2(0,T;\mathbb{R}^2)$, $F\in \mathcal{F}^T$,
$F=\begin{pmatrix}f_1\\f_2\end{pmatrix}$.
\begin{theorem}
The solution to (\ref{wave_eqn})--(\ref{init_eqn}) with a control
$F\in \mathcal{F}^T\cap C_0^\infty(\mathbb{R}_+)$, admits the
following representation:
\begin{equation}
\label{U_1}
u^F(x,t)= \left\{
\begin{array}l\frac{1}{2}f_1(t-x)-\frac{1}{2}f_2(t-x)+\\
+\int_x^t
w_1(x,s)f_1(t-s)+w_2(x,s)f_2(t-s)\,ds,\,0<x<t,\\
-\frac{1}{2}f_1(t+x)-\frac{1}{2}f_2(t+x)+\\
+\int_{-x}^t
w_1(x,s)f_1(t-s)+w_2(x,s)f_2(t-s)\,ds,\,0<-x<t,\\
0, \quad 0<t<|x|.
\end{array}
\right.
\end{equation}
where kernels $w_{1}(x,t)$ and $w_{2}(x,t)$ satisfy the following
Goursat problems:
\begin{eqnarray}
\left\{
\begin{array}l
{w_{1}}_{tt}(x,t)-{w_{1}}_{xx}(x,t)+q(x){w_{1}}(x,t),\quad 0<|x|<t,\\
\frac{d}{dx}{w_{1}}(x,x)=-\frac{q(x)}{4},\quad x>0,\\
\frac{d}{dx}{w_{1}}(x,-x)=-\frac{q(x)}{4},\quad x<0,
\end{array}
\right.\label{W1} \\
\left\{
\begin{array}l
{w_{2}}_{tt}(x,t)-{w_{2}}_{xx}(x,t)+q(x){w_{2}}(x,t),\quad 0<|x|<t,\\
\frac{d}{dx}{w_{2}}(x,x)=\frac{q(x)}{4},\quad x>0,\\
\frac{d}{dx}{w_{2}}(x,-x)=-\frac{q(x)}{4},\quad x<0.
\end{array}
\right. \label{W2}
\end{eqnarray}
\end{theorem}
\begin{proof}
Take arbitrary $F\in \mathcal{F}^T\cap
C_0^\infty(0,T;\mathbb{R}^2)$ and look for $u^F$ in the form
(\ref{U_1}). Then for $x>0$ we have:
\begin{eqnarray*}
u_{xx}(x,t)=\frac{1}{2}f_1''(t-x)-\frac{1}{2}f_2''(t-x)-\frac{d}{dx}w_1(x,x)f_1(t-x)+w_1(x,x)f_1'(t-x)\\
-\frac{d}{dx}w_2(x,x)f_2(t-x)
+w_2(x,x)f_2'(t-x)-{w_1}_x(x,x)f_1(t-x)\\
-{w_2}_x(x,x)f_2(t-x)+\int_{x}^t{w_1}_{xx}(x,s)f_1(t-s)+{w_2}_{xx}(x,s)f_2(t-s)\,ds,\\
u_{tt}(x,t)=\frac{1}{2}f_1''(t-x)-\frac{1}{2}f_2''(t-x)+w_1(x,x)f_1'(t-x)+w_2(x,x)f_2'(t-x)\\
+{w_1}_s(x,x)f_1(t-x)+{w_2}_s(x,x)f_2(t-x)\\
+\int_{x}^t\left({w_1}_{ss}(x,s)f_1(t-s)+{w_2}_{ss}(x,s)f_2(t-s)\right)\,ds,
\end{eqnarray*}
Plugging these expressions into (\ref{wave_eqn}), we obtain that
for $x>0$ the following relation holds true:
\begin{eqnarray}
\label{R_U}0=\int_{x}^t\left(\left({w_1}_{ss}(x,s)-{w_1}_{xx}(x,s)+q(x){w_1}(x,s)\right)f_1(t-s)\right.\\
\left.+\left({w_2}_{ss}(x,s)-{w_2}_{xx}(x,s)+q(x){w_2}(x,s)\right)f_2(t-s)\right)\,ds\notag\\
+f_1(t-x)\left[2\frac{d}{dx}w_1(x,x)+\frac{q(x)}{2}\right]+f_2(t-x)\left[2\frac{d}{dx}w_2(x,x)-\frac{q(x)}{2}\right].\notag
\end{eqnarray}
Similarly, for $x<0$:
\begin{eqnarray*}
u_{xx}(x,t)=-\frac{1}{2}f_1''(t+x)-\frac{1}{2}f_2''(t+x)+\frac{d}{dx}w_1(x,-x)f_1(t+x)+w_1(x,-x)f_1'(t+x)\\
+\frac{d}{dx}w_2(x,-x)f_2(t+x)
+w_2(x,-x)f_2'(t+x)+{w_1}_x(x,-x)f_1(t+x)\\
+{w_2}_x(x,-x)f_2(t+x)+\int_{-x}^t{w_1}_{xx}(x,s)f_1(t-s)+{w_2}_{xx}(x,s)f_2(t-s)\,ds,\\
u_{tt}(x,t)=-\frac{1}{2}f_1''(t+x)-\frac{1}{2}f_2''(t+x)+w_1(x,-x)f_1'(t+x)+w_2(x,-x)f_2'(t+x)\\
+{w_1}_s(x,-x)f_1(t+x)+{w_2}_s(x,-x)f_2(t+x)\\
+\int_{-x}^t\left({w_1}_{ss}(x,s)f_1(t-s)+{w_2}_{ss}(x,s)f_2(t-s)\right)\,ds.
\end{eqnarray*}
Then for $x<0$ we have the equality:
\begin{eqnarray}
\label{L_U}0=\int_{-x}^t\left(\left({w_1}_{ss}(x,s)-{w_1}_{xx}(x,s)+q(x){w_1}(x,s)\right)f_1(t-s)\right.\\
\left.+\left({w_2}_{ss}(x,s)-{w_2}_{xx}(x,s)+q(x){w_2}(x,s)\right)f_2(t-s)\right)\,ds\notag\\
+f_1(t+x)\left[-2\frac{d}{dx}w_1(x,-x)-\frac{q(x)}{2}\right]+f_2(t+x)\left[-2\frac{d}{dx}w_2(x,-x)-\frac{q(x)}{2}\right].\notag
\end{eqnarray}
The condition $\Gamma_0u=F$ at $x=0$ yields that
\begin{eqnarray*}
u^+(\cdot,t)-u^-(\cdot,t)=f_1(t)\\
+\int_0^t\left({w_1}^+(0,s)-{w_1}^-(0,s)\right)f_1(t-s)+\left({w_2}^+(0,s)-{w_2}^-(0,s)\right)f_2(t-s)\,ds,\notag\\
u_x^+(\cdot,t)-u_x^-(\cdot,t)=f_2'(t)\\
+\int_0^t\left({w_1}^+_x(0,s)-{w_1}_x^-(0,s)\right)f_1(t-s)+\left({w_2}^+_x(0,s)-{w_2}^-_x(0,s)\right)f_2(t-s)\,ds,\notag
\end{eqnarray*}
The above equalities imply the continuity of kernels $w_1,$ $w_2$
at $x=0$:
\begin{eqnarray}
{w_1}^+(0,s)-{w_1}^-(0,s),\, {w_2}^+(0,s)-{w_2}^-(0,s),\label{Cont3}\\
{w_1}^+_x(0,s)-{w_1}_x^-(0,s),\,{w_2}^+_x(0,s)-{w_2}_x^-(0,s).\label{Cont4}
\end{eqnarray}
Using the arbitrariness of $F\in \mathcal{F}^T\cap
C_0^\infty(0,T;\mathbb{R}^2)$ in (\ref{R_U}), (\ref{L_U}) and
continuity conditions (\ref{Cont3}), (\ref{Cont3}), we obtain that
$w_1,$ $w_2$ satisfy (\ref{W1}), (\ref{W2}).
\end{proof}
\begin{remark}
When $F\in \mathcal{F}^T$, the function $u^F$ defined by
(\ref{U_1}) is a generalized solution to
(\ref{wave_eqn})--(\ref{init_eqn}).
\end{remark}

The \emph{response operator} $R^T: \mathcal{F}^T\mapsto
\mathcal{F}^T$ with the domain $D_R=\left\{\mathcal{F}^T\cap
C_0^\infty(0,T;\mathbb{R}^2)\right\}$ is defined by
\begin{equation*}
(R^TF)(t):=\left(\Gamma_1u^F\right)(t),\quad 0<t<T.
\end{equation*}
Representation (\ref{U_1}) implies that the response operator has
a form:
\begin{eqnarray}
\label{Resp_repr}
\left(R^TF\right)(t)=\begin{pmatrix}(R_1F)(t)\\(R_2F)(t)\end{pmatrix}= -\frac{1}{2}\begin{pmatrix}f_1'(t)\\-f_2(t)\end{pmatrix}+R*\begin{pmatrix}f_1\\f_2\end{pmatrix}\\
=\begin{pmatrix}-\frac{1}{2}f_1'(t)+\int_0^t\left({w_1}_x(0,s)f_1(t-s)+{w_2}_x(0,s)f_2(t-s)\right)\,ds\\
\frac{1}{2}f_2(t)-\int_0^t\left({w_1}(0,s)f_1(t-s)+{w_2}(0,s)f_2(t-s)\right)\,ds\end{pmatrix},\notag
\end{eqnarray}
where
\begin{equation*}
R(t):=\begin{pmatrix}r_{11}(t) & r_{12}(t)\\
r_{21}(t) & r_{22}(t)\end{pmatrix}=\begin{pmatrix}{w_1}_x(0,t) & {w_2}_x(0,t)\\
-{w_1}(0,t) & -{w_2}(0,t)\end{pmatrix}
\end{equation*}
is a \emph{response matrix}. We introduce the \emph{inner space},
the space of states of system (\ref{wave_eqn})--(\ref{init_eqn})
as $\mathcal{H}^T:=L_2(-T,T)$. The representation (\ref{U_1})
implies that $u^F(\cdot,T)\in \mathcal{H}^T$.

A \emph{control operator} $W^T: \mathcal{F}^T\mapsto
\mathcal{H}^T$ is defined by the formula $W^TF:=u^F(\cdot,T)$. The
\emph{reachable set} is defined by the rule
\begin{equation*}
U^T:=W^T\mathcal{F}^T=\left\{u^F(\cdot,T)\,\big|\,  F\in
\mathcal{F}^T\right\}.
\end{equation*}
We introduce the notations:
\begin{eqnarray*}
S:=\frac{1}{2}\begin{pmatrix} 1 & -1\\ -1 & -1\end{pmatrix},\quad
J^T:\mathcal{F}^T\mapsto \mathcal{F}^T,\quad
\left(J^TF\right)(t)=F(T-t),
\end{eqnarray*}
and note  that
\begin{equation*}
S=S^*,\, SS=\frac{1}{2}I.
\end{equation*}

It will be convenient for us to associate the outer space
$\mathcal{H}^T=L_2(-T,T)$ with a vector space
$L_2(0,T;\mathbb{R}^2)$ by setting for $a\in L_2(-T,T)$ (we keep
the same notation for a function)
\begin{equation*}
a=\begin{pmatrix}a_1(x) \\ a_2(x)\end{pmatrix}\in
L_2(0,T;\mathbb{R}^2),\quad a_1(x):=a(x),\, a_2(x):=a(-x),\, x\in
(0,T).
\end{equation*}
Thus, bearing in mind this association, we consider the control
operator $W^T$, which maps $\mathcal{F}^T$ to
$\mathcal{H}^T=L_2(0,T;\mathbb{R}^2)$, acting (cf. (\ref{U_1})) by
the rule:
\begin{eqnarray*}
\left(W^TF\right)(x)=\begin{pmatrix}
\frac{1}{2}f_1(T-x)-\frac{1}{2}f_2(T-x)\\
-\frac{1}{2}f_1(T-x)-\frac{1}{2}f_2(T-x)
\end{pmatrix}\\
+\begin{pmatrix} \int_x^T
w_1(x,s)f_1(T-s)+w_2(x,s)f_2(T-s)\,ds\\
\int_{x}^T w_1(-x,s)f_1(T-s)+w_2(-x,s)f_2(T-s)\,ds
\end{pmatrix}.
\end{eqnarray*}
On introducing the operator $W: \mathcal{F}^T\mapsto
\mathcal{H}^T=L_2(0,T;\mathbb{R}^2)$ defined by the formula
\begin{equation*}
(WF)(x)=\begin{pmatrix} \int_x^T
w_1(x,s)f_1(s)+w_2(x,s)f_2(s)\,ds\\
\int_{x}^T w_1(-x,s)f_1(s)+w_2(-x,s)f_2(s)\,ds
\end{pmatrix}
\end{equation*}
and noting that $\mathcal{F}^T=\mathcal{H}^T$, we can without
abusing the notations rewrite $W^T$ in a form:
\begin{equation}
\label{WT_form}
W^TF=S\left(I+2SW\right)J^TF=S\left(I+K\right)J^TF,
\end{equation}
where
\begin{equation}
\label{K_W_relations} K=2SW,\quad (KF)(x)=\begin{pmatrix} \int_x^T
k_{11}(x,s)f_1(s)+k_{12}(x,s)f_2(s)\,ds\\
\int_{x}^T k_{21}(x,s)f_1(s)+k_{22}(x,s)f_2(s)\,ds
\end{pmatrix}.
\end{equation}

\begin{theorem}
\label{ControlTheor} The control operator is a boundedly
invertible isomorphism between $\mathcal{F}^T$ and
$\mathcal{H}^T$, and $U^T=\mathcal{H}^T$.
\end{theorem}
\begin{proof}
It is clear that in representation (\ref{WT_form}) each of the
operators $S: \mathcal{H}^T\mapsto \mathcal{H}^T,$ $I+K:
\mathcal{F}^T\mapsto \mathcal{H}^T$, $J^T: \mathcal{F}^T\mapsto
\mathcal{F}^T$ is boundedly invertible isomorphism.
\end{proof}

The \emph{connecting operator} $C^T:\mathcal{F}^T\mapsto
\mathcal{F}^T$ is introduced via the quadratic form:
\begin{equation*}
\left(C^T
F_1,F_2\right)_{\mathcal{F}^T}=\left(u^{F_1}(\cdot,T),u^{F_2}(\cdot,T)\right)_{\mathcal{H}^T}.
\end{equation*}
The crucial fact in the Boundary Control method is that the
connecting operator is expressed in terms of inverse dynamic data:
\begin{theorem}
The connecting operator $C^T$ admits the following representation:
\begin{equation*}
\left(C^TF\right)(t)=\frac{1}{2}\begin{pmatrix}f_1(t)\\
f_2(t)\end{pmatrix}+\int_0^TC(t,s)\begin{pmatrix}
f_1(s)\\f_2(s)\end{pmatrix}\,ds,
\end{equation*}
where
\begin{align*}
&C_{1,1}(t,s)= p_1(2T-t-s)-p_1(|t-s|),\quad p_1(s)=\int_0^s
r_{11}(\alpha)\,d\alpha,\\
&C_{1,2}(t,s)=\widetilde p_1(2T-t-s)-\widetilde p_1(t-s),\quad
\widetilde p_1(s)=\left\{\begin{array}l \int_0^s
r_{12}(\alpha)\,d\alpha,\, s>0,\\
-\int_0^{-s} r_{12}(\alpha)\,d\alpha,\, s<0,
\end{array}
\right.\\
&C_{2,1}(t,s)=-\widetilde r_{21}(t-s)-\widetilde
r_{21}(2T-t-s),\quad \widetilde r_{21}(s)=\left\{\begin{array}l
r_{21}(s),\, s>0,\\
-r_{21}(-s),\, s<0,
\end{array}
\right.\\
&C_{2,2}(t,s)=-r_{22}(|t-s|)-r_{22}(2T-t-s).
\end{align*}
\end{theorem}
\begin{proof}
We take $F,G\in \mathcal{F}^T\cap C_0^\infty(0,T;\mathbb{R}^2)$
and introduce the Blagoveschenskii function by setting
\begin{equation*}
\Psi(t,s)=\left(u^F(\cdot,t),u^G(\cdot,s)\right)_{\mathcal{H}^T},\quad
s,t>0.
\end{equation*}
We show that $\Psi$ satisfy the wave equation. Indeed, using that
$u^F_{tt}=-Hu^F$ and Green identity, we can evaluate:
\begin{eqnarray*}
\Psi_{tt}(t,s)-\Psi_{ss}(t,s)=\left(-H^*u^F(\cdot,t),u^G(\cdot,s)\right)_{\mathcal{H}^T}+\left(u^F(\cdot,t),H^*u^G(\cdot,s)\right)_{\mathcal{H}^T}\\
=\left(\left(\Gamma_0u^F\right)(t),\left(\Gamma_1u^G\right)(s)
\right)_B-\left(\left(\Gamma_1u^F\right)(t),\left(\Gamma_0u^G\right)(s)
\right)_B=:P(t,s).
\end{eqnarray*}
Note that $\Psi$ satisfy $\Psi(0,s)=\Psi_t(0,s)=0,$ and that
\begin{equation*}
\Psi(T,T)=\left(u^F(\cdot,T),u^G(\cdot,T)\right)_{\mathcal{H}^T}=\left(C^TF,G\right)_{\mathcal{F}^T}.
\end{equation*}
So, by d'Alembert formula:
\begin{equation}
\label{CT_integr}
\left(C^TF,G\right)_{\mathcal{F}^T}=\int_0^T\int_{\tau}^{2T-\tau}P(\tau,\sigma)\,d\sigma\,d\tau.
\end{equation}
We rewrite the right hand side:
\begin{equation}
\label{P_RHS}
P(t,s)=\left(\begin{pmatrix}f_1(t)\\
f_2'(t)\end{pmatrix},(RG)(s)\right)_{B}-\left((RF)(t),\begin{pmatrix}g_1(s)\\
g_2'(s)\end{pmatrix}\right)_{B},
\end{equation}
and continue the functions $g_1,$ $g_2$ (we keep the same
notations) from $(0,T)$ to the interval $(0,2T)$ by the rule:
\begin{equation}
\label{G_form}
g_1(s)=\left\{\begin{array}l g_1(s),\, 0<s<T,\\
-g_1(2T-s),\, T<s<2T,
\end{array}
\right.
g_2(s)=\left\{\begin{array}l g_2(s),\, 0<s<T,\\
g_2(2T-s),\, T<s<2T.
\end{array}
\right.
\end{equation}
After such a continuation the second term in (\ref{P_RHS}) become
odd in $s$ with respect to $s=T$ and disappears after integration
in (\ref{CT_integr}), so we come to the following expression
\begin{equation}
\label{CT_1}
\left(C^TF,G\right)_{\mathcal{F}^T}=\int_0^T\int_{\tau}^{2T-\tau}\left(\begin{pmatrix}f_1(t)\\
f_2'(t)\end{pmatrix},(RG)(s)\right)_{B}\,d\sigma\,d\tau.
\end{equation}
Integrating by parts in (\ref{CT_1}) and using that
$C^T=\left(C^T\right)^*$ and arbitrariness of $F$ yields
\begin{equation}
\label{CT_2}
\left(C^TG\right)(\tau)=\begin{pmatrix}\int_{\tau}^{2T-\tau}(R_1G)(\sigma)\,d\sigma\\
(R_2G)(\tau)+(R_2G)(2T-\tau)\end{pmatrix}.
\end{equation}
Evaluating (\ref{CT_2}) making use of (\ref{Resp_repr}) and
continuation of $G$ (\ref{G_form}),  we obtain that
\begin{eqnarray}
\left(C^TG\right)(\tau)=\frac{1}{2}\begin{pmatrix}g_1(\tau)\\
g_2(\tau)\end{pmatrix}
+\frac{1}{2}\begin{pmatrix}\int_{\tau}^{2T-\tau}\int_0^\sigma
\left(r_{11}(s)g_1(\sigma-s)+r_{12}(s)g_2(\sigma-s)\right)\,d\sigma\\
-\int_0^\tau\left(r_{21}(s)g_1(\tau-s)+r_{22}(s)g_2(\tau-s)\right)\,ds\end{pmatrix}\notag\\
+\begin{pmatrix}0\\\int_0^{2T-\tau}\left(r_{21}(s)g_1(2T-\tau-s)+r_{22}(s)g_2(2T-\tau-s)\right)\,ds\end{pmatrix}.\label{CT_3}
\end{eqnarray}
Consider the term
\begin{equation}
\label{C_11} \int_{\tau}^{2T-\tau}\int_0^\sigma
r_{11}(s)g_1(\sigma-s)\,ds\,d\sigma=I(2T-\tau)-I(\tau),
\end{equation}
where
\begin{equation*}
I(\tau)=\int_0^\tau\int_\alpha^\tau
r_{11}(\sigma-\alpha)g_1(\alpha)\,d\sigma\,d\alpha.
\end{equation*}
We evaluate (\ref{C_11}) using that $g_1$ is odd with respect to
$T$:
\begin{equation}
\label{I1} I(\tau)=\int_0^\tau\int_0^{|\tau-\alpha|}
r_{11}(\sigma)\,d\sigma g_1(\alpha)\,d\alpha=\int_0^\tau
p_1(|\tau-\alpha|)g_1(\alpha)\,d\alpha,
\end{equation}
where $p_1(s)=\int_0^s r_{11}(\alpha)\,d\alpha.$ The first term in
(\ref{C_11}) we can rewrite in a form:
\begin{eqnarray}
I(2T-\tau)=\left(\int_0^T+\int_\tau^{2T-\tau}\right)\int_0^{2T-\tau-\alpha}
r_{11}(\sigma)\,d\sigma g_1(\alpha)\,d\alpha\notag\\
=\int_0^T p_1(2T-\tau-\alpha)g_1(\alpha)\,d\alpha-\int_{\tau}^T
p_1(\alpha-\tau)g_1(\alpha)\,d\alpha. \label{I2}
\end{eqnarray}
Then from (\ref{I1}) and (\ref{I2}) we obtain that
\begin{equation*}
\int_{\tau}^{2T-\tau}\int_0^\sigma
r_{11}(s)g_1(\sigma-s)\,ds\,d\sigma=\int_0^T
\left(p_1(2T-\tau-\alpha)-p_1(|\alpha-\tau|)g_1(\alpha)\right)\,d\alpha,
\end{equation*}
which proves the formula for $C_{11}$. Now we consider the term
\begin{equation}
\label{C_12} \int_{\tau}^{2T-\tau}\int_0^\sigma
r_{12}(s)g_2(\sigma-s)\,ds\,d\sigma.
\end{equation}
Note that it has the same structure as (\ref{C_11}), but we should
take into account that $g_2$ is odd with respect to $T$. Counting
this, we have that:
\begin{equation*}
\label{I3} I(2T-\tau)=\int_0^T
p_2(2T-\tau-\alpha)g_2(\alpha)\,d\alpha+\int_{\tau}^T
p_2(\alpha-\tau)g_2(\alpha)\,d\alpha,
\end{equation*}
where $p_2(s)=\int_0^s r_{12}(\alpha)\,d\alpha.$ Then
\begin{eqnarray}
I(2T-\tau)-I(\tau)=\int_0^T
p_2(2T-\tau-\alpha)g_2(\alpha)\,d\alpha\label{I4}\\
+\int_{\tau}^T p_2(\alpha-\tau)g_2(\alpha)\,d\alpha-\int_{0}^T
p_2(|\alpha-\tau|)g_2(\alpha)\,d\alpha,\notag
\end{eqnarray}
After we introduce the notation
\begin{equation*}
\widetilde p_1(s)=\left\{\begin{array}l \int_0^s
r_{12}(\alpha)\,d\alpha,\, s>0,\\
-\int_0^{-s} r_{12}(\alpha)\,d\alpha,\, s<0,
\end{array}
\right.=\left\{\begin{array}l p_2(s),\, s>0,\\
-p_2(-s),\, s<0,
\end{array}
\right.
\end{equation*}
we can rewrite (\ref{C_12}), taking into account (\ref{I4}), as
\begin{equation*}
\int_{\tau}^{2T-\tau}\int_0^\sigma
r_{12}(s)g_2(\sigma-s)\,ds\,d\sigma=\int_0^T \left(\widetilde
p_1(2T-\tau-\alpha)-\widetilde
p_1(\tau-\alpha)\right)g_2(\alpha)\,d\alpha,
\end{equation*}
which proves the formula for $C_{12}$. Similarly one can prove
formulaes for $C_{21},$ $C_{22}$.

\end{proof}

We note that the symmetry of $C^T$ implies the restriction on the
entries, specifically, the following relation holds:
\begin{equation*}
C_{2,1}(t,s)=C_{1,2}(t,s).
\end{equation*}
This equality is equivalent to
\begin{equation*}
-\widetilde r_{21}(t-s)-\widetilde r_{21}(2T-t-s)=\widetilde
p_1(2T-t-s)-\widetilde p_1(s-t),
\end{equation*}
which yields:
\begin{equation*}
-\widetilde r_{21}(s)=\widetilde p_1(s).
\end{equation*}

\begin{remark}
The components of the response matrix have to be connected by the
relation:
\begin{equation*}
r_{21}'(s)=-r_{12}(s),\quad s>0.
\end{equation*}
\end{remark}

\section{Dynamic inverse problem}

In this section we derive equations of inverse dynamic problem,
using them we answer the question on recovering a potential
$q(x)$, $x\in (-T,T)$ from the response operator $R^{2T}$.

\subsection{Krein equations.}

Let $y(x)$ be a solution to the following Cauchy problem:
\begin{equation}
\label{Cauchy_pr} \left\{
\begin{array}l
-y''+qy=0,\quad x\in (-T,T),\\
y(0)=0,\,y'(0)=1.
\end{array}
\right.
\end{equation}

We set up the \emph{special control problem}: to find $F\in
\mathcal{F}^T$ such that $W^TF=y$ in $\mathcal{H}^T$. By the
Theorem \ref{ControlTheor}, such a control $F$ exists, but we can
say even more:
\begin{theorem}
The solution to a special control problem is a unique solution to
the following equation:
\begin{equation}
\label{Krein_eqn}
\left(C^TF\right)(t)=(T-t)\begin{pmatrix} 1\\
0\end{pmatrix},\quad t\in (0,T).
\end{equation}
\end{theorem}
\begin{proof}
We observe that if $G\in \mathcal{F}^T\cap
C_0^\infty(0,T;\mathbb{R}^2)$, then it is true that
\begin{equation*}
u^G(x,T)=\int_0^T(T-t)u^G_{tt}(x,t)\,dt.
\end{equation*}
Using this, we can evaluate the quadratic form:
\begin{eqnarray*}
\left(C^TF,G\right)_{\mathcal{F}^T}=\left(W^TF,W^TG\right)_{\mathcal{H}^T}=\left(y(\cdot),u^G(\cdot,T)\right)_{\mathcal{H}^T}\\
=\int_{-T}^Ty(x)\int_0^T(T-t)u^G_{tt}(x,t)\,dt\,dx=\int_{0}^T(T-t)\left(y(\cdot),
-H^*u^G(\cdot,t)\right)_{\mathcal{H}^T}\,dx\,dt\\
=\int_{0}^T(t-T)\left[\left(\left(\Gamma_0y(\cdot)\right)(t),\left(\Gamma_1u^G\right)(t)\right)_B-\left(\left(\Gamma_1y(\cdot)\right)(t),\left(\Gamma_0u^G\right)(t)\right)_B\right]\,dt\\
=\int_0^T(T-t)\left(\begin{pmatrix}
1\\0\end{pmatrix},\begin{pmatrix}g_1(t)\\g_2'(t)\end{pmatrix}\right)\,dt,
\end{eqnarray*}
from where (\ref{Krein_eqn}) follows due to the arbitrariness of
$G$.
\end{proof}
Representation formulas (\ref{U_1})  imply that  that the solution
$F$ to a special control problem satisfies relations:
\begin{eqnarray*}
y(T)=u^F(T,T)=\frac{1}{2}f_1(0)-\frac{1}{2}f_2(0),\\
y(-T)=u^F(-T,T)=-\frac{1}{2}f_1(0)-\frac{1}{2}f_2(0).
\end{eqnarray*}
Thus solving (\ref{Krein_eqn}) for all $T\in (0,T)$, we recover
the solution $y(x)$ to (\ref{Cauchy_pr}) on the interval $(-T,T)$.
Then the potential $q(x),$ $x\in (-T,T)$ can be recovered as
$q(x)=\frac{y''(x)}{y(x)}$, $x\in (-T,T)$.

\subsection{Gelfand-Levitan equations}

We introduce the notation:
\begin{eqnarray}
C^T=\frac{1}{2}(I+C),\quad
(Cf)(t)=2\int_0^TC(t,s)\label{C_def}\begin{pmatrix}
f(s)\\G(s)\end{pmatrix}\,ds.
\end{eqnarray}
For $F,G\in \mathcal{F}^T$ we set $W^TF=a$, $W^TG=b$, where
$a,b\in \mathcal{H}^T$, on using the controllability (Theorem
\ref{ControlTheor}), we have that (see (\ref{WT_form}))
\begin{eqnarray*}
F=J^T(I+K)^{-1}S^{-1}a=2J^T(I+K)^{-1}Sa,\\
G=J^T(I+K)^{-1}S^{-1}b=2J^T(I+K)^{-1}Sb.
\end{eqnarray*}
Using above representations we can rewrite the quadratic form as:
\begin{eqnarray}
\left(C^TF,G\right)_{\mathcal{H}^T}=\left(\frac{1}{2}(I+C)2J^T(I+K)^{-1}Sa,2J^T(I+K)^{-1}Sb\right)_{\mathcal{H}^T}\notag\\
=\left(2\left((I+K)^{-1}\right)^*J^T(I+C)J^T(I+K)^{-1}Sa,Sb\right)_{\mathcal{H}^T}.\label{C1}
\end{eqnarray}
On the other hand:
\begin{equation}
\label{C2}
\left(C^TF,G\right)_{\mathcal{H}^T}=\left(W^TF,W^TG\right)_{\mathcal{H}^T}=(a,b)_{\mathcal{H}^T}=(2Sa,Sb)_{\mathcal{H}^T}.
\end{equation}
On comparing (\ref{C1}) and (\ref{C2}), we obtain the following
operator identity:
\begin{equation}
\label{OperId} \left((I+K)^{-1}\right)^*J^T(I+C)J^T(I+K)^{-1}=I.
\end{equation}
We introduce the following notations
\begin{eqnarray}
\label{M_oper}I+M=(I+K)^{-1},\\
\left(MF\right)(x)=\begin{pmatrix} \int_x^T
m_{11}(x,s)f_1(s)+m_{12}(x,s)f_2(s)\,ds\notag \\
\int_{x}^T m_{21}(x,s)f_1(s)+m_{22}(x,s)f_2(s)\,ds
\end{pmatrix}\\
\left(M^*a\right)(t)= \begin{pmatrix} \int_0^t
m_{11}(x,t)a_1(x)+m_{21}(x,t)a_2(s)\,dx\\
\int_{t}^t m_{12}(x,t)a_1(s)+m_{22}(x,t)a_2(x)\,dx
\end{pmatrix}.\notag
\end{eqnarray}
It is easy to check that on a diagonal the kernels of $K,$ $M$
satisfy a relation
\begin{equation}
\label{M_K_relations} m_{ij}(x,x)=-k_{ij}(x,x),\quad
i,j=\{1,2\},\quad x\in (0,T).
\end{equation}
In new notations the operator equality (\ref{OperId}) has a form:
\begin{equation}
\label{OperId0} (I+M)^*(I+\widetilde C)(I+M)=I,
\end{equation}
where
\begin{equation}
\label{C_wid_def} \widetilde C=J^TCJ^T,\quad \left(\widetilde
CF\right)(t)=\int_0^T \widetilde C(t,s)F(s)\,ds.
\end{equation}
We rewrite (\ref{OperId0}) in a form:
\begin{equation}
\label{OperId1} M^*+(I+M)^*\left(M+\widetilde C+\widetilde C
M\right)=0
\end{equation}
On introducing a function
\begin{equation*}
\Phi(x,s)=m(x,s)+\widetilde C(x,s)+\int_0^T\widetilde
C(x,\alpha)m(\alpha,s)\,d\alpha,\quad x,s\in (0,T),
\end{equation*}
we can write down an equality on the kernel for the operator in
the left hand side in (\ref{OperId1}) $M^*+\Phi+M^*\Phi=0$:
\begin{equation*}
m(s,x)+\Phi(x,s)+\int_0^tm(\alpha,x)\Phi(\alpha,s)\,d\alpha=0,\quad
x,s\in(0,T).
\end{equation*}
Since $m(s,x)=0$ when $x<s$, we obtain that $\Phi$ satisfies the
relation:
\begin{equation*}
\Phi(x,s)+\int_0^tm(\alpha,x)\Phi(\alpha,s)\,d\alpha=0,\quad x<s.
\end{equation*}
Due to the fact that $\Psi$ satisfies  a Volterra equation of a
second kind, we obtain that $\Phi(x,s)=0$ for $x<s$, which
immediately yields the following equation on the matrix function
$m$:
\begin{equation}
\label{G_L_eqn} m(x,s)+\widetilde C(x,s)+\int_0^T\widetilde
C(x,\alpha)m(\alpha,s)\,d\alpha=0,\quad 0<x<s<T.
\end{equation}
As a result we can formulate the following
\begin{theorem}
The matrix kernel of the operator $M$ (\ref{M_oper}) satisfy the
Gelfand-Levitan equation (\ref{G_L_eqn}), where the kernel
$\widetilde C$ is defined by (\ref{C_def}), (\ref{C_wid_def}) by
solving which, one can recover the potential using relations
between kernels (\ref{K_W_relations}), (\ref{M_K_relations}) and
relations on diagonals $\{x=t\},$ $\{-x=t\}$ in (\ref{W1}),
(\ref{W2}):
\begin{eqnarray*}
q(x)=2\frac{d}{dx}\left(m_{11}(x,x)-m_{12}(x,x)\right),\quad x\in (0,T),\\
q(-x)=-2\frac{d}{dx}\left(m_{11}(x,x)+m_{12}(x,x)\right),\quad
x\in (0,T).
\end{eqnarray*}
\end{theorem}

\subsection{Relationship between dynamic and spectral inverse data. }

The problem of finding relationships between different types of
inverse data is very important in inverse problems theory. We can
mention \cite{B01JII,B03,AMR,MM3,MM4} on some recent results in
this direction. Below we show the relationship between the dynamic
response function and matrix spectral measure.

Consider two solution to the equation
\begin{equation}
\label{Eq_Sch} -\phi''+q(x)\phi=\lambda\phi,\quad
-\infty<x<\infty,
\end{equation}
satisfying the Cauchy data:
\begin{equation*}
\varphi(0,\lambda)=0,\, \varphi'(0,\lambda)=1,\,
\theta(0,\lambda)=-1,\,\theta'(0,\lambda)=0.
\end{equation*}
Note that
\begin{equation*}
\Gamma_0\varphi=0,\,\Gamma_0\theta=0,\,
\Gamma_1\varphi=\begin{pmatrix}1\\ 0\end{pmatrix},\,
\Gamma_1\theta=\begin{pmatrix}0\\ 1\end{pmatrix}.
\end{equation*}
We fix some $N>0$ and prescribe self-adjoint boundary conditions
at $x=\pm N:$
\begin{eqnarray}
a_1\phi(-N,\lambda)+b_1\phi'(-N,\lambda)=0,\label{Eq_Sch1} \\
a_2\phi(N,\lambda)+b_2\phi'(N,\lambda)=0.\label{Eq_Sch2}
\end{eqnarray}
Eigenvalues and normalized eigenfunctions of (\ref{Eq_Sch}),
(\ref{Eq_Sch1}), (\ref{Eq_Sch2}) are denoted by
$\{\lambda_n,y_n\}_{n=1}^\infty.$ Let $\beta_n,\gamma_n\in
\mathbb{R}$ be such that
\begin{equation*}
y_n(x)=\beta_n\varphi(x,\lambda_n)+\gamma_n\theta(x,\lambda_n),
\quad \text{then}\quad \Gamma_1y_n=\begin{pmatrix}\beta_n\\
\gamma_n\end{pmatrix}.
\end{equation*}
Let $F\in \mathcal{F}^T\cap C^\infty_0(0,T;\mathbb{R}^2)$, and
$v^F$ be a solution to (\ref{wave_eqn})--(\ref{init_eqn}),
(\ref{Eq_Sch1}), (\ref{Eq_Sch2}), i.e. a solution to the wave
equation on the interval $(-N,N)$. On multiplying wave equation
for $v^F$ by $y_n$ and integrating by parts, we get the following
relation:
\begin{eqnarray*}
0=\int_{-T}^T
v^F_{tt}y_n\,dx-\int_{-N}^Nv^F_{xx}y_n\,dx+\int_{-N}^N
q(x)v^Fy_n\,dx=\int_{-N}^N v^F_{tt}y_n\,dx\\
+\left(v^F,Hy_n\right)+\left(\Gamma_1v^F,\Gamma_0y_n\right)_B-\left(\Gamma_0v^F,\Gamma_1y_n\right)_B\\
=\int_{-T}^T
v^F_{tt}y_n\,dx+\lambda_n\left(v^F,y_n\right)-\left(\begin{pmatrix}f_1(t)\\ f_2'(t)\end{pmatrix},\begin{pmatrix}\beta_n\\
\gamma_n\end{pmatrix}\right)_B.
\end{eqnarray*}
Looking for the solution to (\ref{wave_eqn})--(\ref{init_eqn}) in
a form
\begin{equation}
\label{spectr_repr} v^F=\sum_{k=1}^\infty c_k(t)y_k(x),
\end{equation}
we plug (\ref{spectr_repr}) into (\ref{wave_eqn}) and multiply by
$y_n$ to get:
\begin{equation*}
\int_{-N}^N
\sum_{k=1}^\infty c_k''(t)y_k(x)y_n(x)\,dx+\int_{-N}^N\sum_{k=1}^\infty c_k(t)y_k(x)\lambda_ny_n(x)\,dx=\left(\begin{pmatrix}f_1(t)\\ f_2'(t)\end{pmatrix},\begin{pmatrix}\beta_n\\
\gamma_n\end{pmatrix}\right)_B.
\end{equation*}
Thu we obtain that $c_n(t)$, $n\geqslant 1,$ satisfies the
following Cauchy problem:
\begin{equation*}
\left\{ \begin{array}l c_n''(t)+\lambda_nc_n(t)=\left(\begin{pmatrix}f_1(t)\\ f_2'(t)\end{pmatrix},\begin{pmatrix}\beta_n\\
\gamma_n\end{pmatrix}\right)_B,\\
c_n(0)=0,\, c_n'(0)=0.
\end{array}
\right.
\end{equation*}
the solution of which is given by the formula
\begin{equation*}
c_n(t)=\int_0^t\frac{\sin{\sqrt{\lambda_n}(t-s)}}{\sqrt{\lambda_n}}\left(f_1(s)\beta_n+f_2'(s)\gamma_n\right)\,ds.
\end{equation*}
Then for $v^F$ (\ref{spectr_repr}) we have the expansion:
\begin{eqnarray}
v^F(x,t)=\sum_{k=1}^\infty\int_0^t\frac{\sin{\sqrt{\lambda_n}(t-s)}}{\sqrt{\lambda_n}}\left(f_1(s)\beta_n+f_2'\gamma_n\right)\,ds\left(\beta_n\varphi(x,\lambda_n)+\gamma_n\theta(x,\lambda_n)\right)\notag\\
=\sum_{k=1}^\infty\int_0^t\frac{\sin{\sqrt{\lambda_n}(t-s)}}{\sqrt{\lambda_n}}\left(\begin{pmatrix}\beta_n\\
\gamma_n\end{pmatrix}\otimes \begin{pmatrix}\beta_n\\
\gamma_n\end{pmatrix}\begin{pmatrix}f_1(s)
\\ f_2'(s)\end{pmatrix},\begin{pmatrix}\varphi(x,\lambda_n)\\
\theta(x,\lambda_n)\end{pmatrix}\right)\notag\\
=\int_{-\infty}^\infty\int_0^t\frac{\sin{\sqrt{\lambda}(t-s)}}{\sqrt{\lambda}}\left(d\Sigma_N(\lambda)\begin{pmatrix}f_1(s)
\\ f_2'(s)\end{pmatrix},\begin{pmatrix}\varphi(x,\lambda)\\
\theta(x,\lambda)\end{pmatrix}\right).\label{U_F_repr}
\end{eqnarray}
Where $d\Sigma_N(\lambda)$ is a matrix measure (see \cite{Le}).
Due to the finite speed of the wave propagation in system
(\ref{wave_eqn})--(\ref{init_eqn}) (equal to one), we have that
\begin{equation}
\label{Func_eq} v^F(\cdot,t)=u^F(\cdot,t),\quad \text{for}\, t<N,
\end{equation}
and for $T<N$ holds that $R^{2T}F=\Gamma_1v^F$. Thus the response
operator $R^T$ for $T<2N$, is given by
\begin{eqnarray}
(RF)(t)=\Gamma_1v^F=\sum_{k=1}^\infty c_k(t)\Gamma_1y_k=\sum
c_k(t)\begin{pmatrix} \beta_k\\
\gamma_k\end{pmatrix}\label{R_T_reprSP}\\
=\sum_{k=1}^\infty\int_0^t\frac{\sin{\sqrt{\lambda_k}(t-s)}}{\sqrt{\lambda_k}}\left(f_1(s)\beta_k+f_2'\gamma_k\right)\,ds\begin{pmatrix} \beta_k\\
\gamma_k\end{pmatrix}\notag\\
=\int_{-\infty}^\infty\int_0^t
\frac{\sin{\sqrt{\lambda}(t-s)}}{\sqrt{\lambda}}d\Sigma_N(\lambda)\begin{pmatrix}
f_1(s)\\ f_2'(s)\end{pmatrix}\,ds,\quad 0<t<2N.\notag
\end{eqnarray}

Taking $F,G\in \mathcal{F}^T\cap C^\infty_0(0,T;\mathbb{R}^2)$,
for $T<N$ we evaluate the form using (\ref{U_F_repr}) and
(\ref{Func_eq}):
\begin{eqnarray}
(C^TF,G)_{\mathcal{F}^T}=(u^F,u^G)_{\mathcal{H}^T}=(v^F,v^G)_{\mathcal{H}^T}\label{C_T_reprSP}\\
=\sum_{k=1}^\infty\int_0^T\int_{0}^T\frac{\sin{\sqrt{\lambda_n}(t-s)}}{\sqrt{\lambda_n}}\left(f_1\beta_n+f_2'\gamma_n\right)\,ds
\frac{\sin{\sqrt{\lambda_n}(t-\tau)}}{\sqrt{\lambda_n}}\left(g_1\beta_n+g_2'\gamma_n\right)\,d\tau
\notag   \\
=\int_0^T\int_{0}^T\int_{-\infty}^\infty\frac{\sin{\sqrt{\lambda}(t-s)}}{\sqrt{\lambda}}
\frac{\sin{\sqrt{\lambda}(t-\tau)}}{\sqrt{\lambda}}\left(d\Sigma_N(\lambda)\begin{pmatrix}f_1(s)
\\ f_2'(s)\end{pmatrix},\begin{pmatrix}g_1(\tau)
\\ g_2'(\tau)\end{pmatrix}\right)\,ds\,d\tau\notag
\end{eqnarray}

We observe that due to the finite speed of wave propagation in
system (\ref{wave_eqn})--(\ref{init_eqn}) equal to one, in
representation formulas for response operator (\ref{R_T_reprSP})
and for connecting operator (\ref{C_T_reprSP}) we can substitute
$d\Sigma_N(\lambda)$ by any $d\Sigma_M(\lambda)$, $M>N$, where
$d\Sigma_M(\lambda)$ corresponds to some selfadjoint boundary
conditions at $\pm M$, or we can let $N$ go to infinity, and
substitute by a limit measure $d\Sigma(\lambda)$ (see \cite{Le}).

The inverse problem for a Schr\"odinger operator on a half-line
from a spectral measure is solved in \cite{GL}, in \cite{Le} the
inverse spectral problem for a Schr\"odinger operator on a real
line is discussed, but some questions remain open. The authors in
\cite{AM,AMR,MM3} point out the relationships between the dynamic
and spectral inverse problems for the case of half-line.

\begin{remark}
The control, response and connecting operators admit
representations in terms of the spectral inverse data (matrix
measure $d\Sigma$), see (\ref{U_F_repr}), (\ref{R_T_reprSP}) and
(\ref{C_T_reprSP}). This circumstance makes it possible to apply
the Boundary Control method in the spirit of \cite{AM,AMR,MM3} for
studying the inverse spectral problem from a matrix measure.
\end{remark}

\noindent{\bf Acknowledgments}

The research of Victor Mikhaylov was supported in part by RFBR
17-01-00529. Alexandr Mikhaylov was supported by RFBR 17-01-00099;
A. S. Mikhaylov and V. S. Mikhaylov were partly supported by VW
Foundation program "Modeling, Analysis, and Approximation Theory
toward application in tomography and inverse problems."

\end{document}